\documentclass{article}
\usepackage{amssymb,amsmath,amsthm}
\newtheorem{theorem}{Theorem}
\newtheorem{definition}{Definition}
\newtheorem{lemma}{Lemma}

\newtheorem{remark}{Remark}
\date{}
\numberwithin{equation}{section} \numberwithin{theorem}{section}
\numberwithin{lemma}{section} \numberwithin{corollary}{section}
\numberwithin{remark}{section} \numberwithin{proposition}{section}
\numberwithin{definition}{section}
\def\avint{\mathop{\,\rlap{--}\!\!\int}\nolimits}
\begin{document}
\newcommand{\n}{\noindent}
\newcommand{\vs}{\vskip}
\title{Lipschitz Continuity for
Elliptic Free Boundary Problems with Dini Mean Oscillation Coefficients}

\vskip 0.5cm
\author{A. Lyaghfouri \\
United Arab Emirates University\\
Department of Mathematical Sciences\\
Al Ain, Abu Dhabi, UAE} \maketitle

 \vskip 0.5cm
\begin{abstract}
We establish local interior Lipschitz continuity of the solutions of a class of
free boundary elliptic problems assuming the coefficients of the equation of
Dini mean oscillation in at least one direction. The novelty in this regularity 
result lies in the fact that it allows discontinuous coefficients in all but one variable.
\end{abstract}

\vs 0.5cm

\n Key words: Dini Mean Oscillation Condition, Free boundary Problem, Harnack's inequality, Lipschitz Continuity.

\vs 0.5cm \noindent AMS Mathematics Subject Classification:
35J15, 35R35
\section{Introduction}\label{1}

\vskip 0.3cm \n Throughout this paper, we denote by $\Omega$ a bounded domain in $\mathbb{R}^n$ and by
$\mathbf{A}(x)= (a_{ij}(x))$ an $n\times n$ matrix
that satisfies for some positive constant $\lambda\in(0,1)$
\begin{eqnarray}\label{e1.1-1.2}
&\sum_{i,j}|a_{ij}(x)|\leq \lambda^{-1},\quad\text{for a.e. }x\in \Omega,\\
&\mathbf{A}(x)\cdot\xi\cdot\xi \geq \lambda\vert \xi\vert^2,\quad\text{for a.e. }x\in \Omega,
\quad\text{for all }\xi\in\mathbb{R}^n
\end{eqnarray}

\n $\mathbf{f}: \Omega \rightarrow \mathbb{R}^n$ is a vector function such that 
$\mathbf{f}(x)=(f_1(x),...,f_n(x))$ and
\begin{equation}\label{e1.3}
\exists \overline{f}>0:~~ ||\mathbf{f}||_\infty\leq \overline{f}
\end{equation}

\n We consider the following problem

\begin{equation*}(P)
\begin{cases}
& \text{Find } (u, \chi) \in  H^{1}(\Omega)\times L^\infty (\Omega)
\text{ such that}:\\
& (i)\quad  u\geq 0, \quad 0\leq \chi\leq 1 ,
\quad u(1-\chi) = 0 \,\,\text{ a.e.  in } \Omega\\
& (ii)\quad \text{div}(\mathbf{A}(x) \nabla
u + \chi \mathbf{f}(x))=0\quad \text{ in } H_0^{-1}(\Omega)\\
\end{cases}
\end{equation*}

\n This class covers a set of various problems including the heterogeneous dam problem \cite{[A1]} \cite{[ChiL1]}, \cite{[CL2]},
\cite{[L1]}, in which case $\Omega$ represents a porous medium with permeability matrix $\mathbf{A}(x)$,
and $\mathbf{f}(x)=\mathbf{A}(x)\mathbf{e}$, with $\mathbf{e}=(0,...,0,1)$.
A second example is the lubrication problem \cite{[AC]} which is
obtained when $\mathbf{A}(x)=h^3(x)\mathbf{I_2}$ and $\mathbf{f}(x)=h(x)\mathbf{e}$, where $\mathbf{I_2}$ is the
$2\times2$ identity matrix,
and $h(x)$ is a scalar function related to the Reynolds equation.
A third example is the aluminium electrolysis problem \cite{[BMQ]} which corresponds to
$\mathbf{A}(x)=k(x)\mathbf{I_2}$ and $\mathbf{f}(x)=h(x)\mathbf{e}$, with $k(x)$ and
$h(x)$ two given scalar functions.

\vs 0.3cm\n We observe that if $\mathbf{f}\in L_{loc}^q(\Omega)$ for some $q>n$, then so is $\chi \mathbf{f}$, and by
taking into account the assumptions (1.1)-(1.2) and the equation $(P)ii)$, we infer from \cite{[GT]} Theorem 8.24, p. 202
that $u\in C_{loc}^{0,\alpha}( \Omega)$ for any $\alpha\in(0,1)$.
In this paper, we will improve this regularity by showing that under suitable assumptions, we actually have $u\in C_{loc}^{0,1}(\Omega)$.
We observe that this regularity is optimal due to the gradient discontinuity across the free boundary
which is the interface that separates the sets $\{u=0\}$ and $\{u>0\}$ from each other. 
Moreover, Lipschitz continuity is not only interesting by itself, but is also of particular importance in the 
analysis of the free boundary (see for example \cite{[CL1]} and \cite{[CL3]}).

\vs 0,3cm \n Before stating our main result, we need to introduce a definition.

\begin{definition}\label{d2.1}
i)~ We say that a function $\omega: (0,1] \rightarrow [0,\infty)$ satisfies the Dini condition if 
\[\int_0^1 \frac{\omega(r)}{r}dr <\infty\]

\n ii)~ A function $f\in L^1(\Omega)$ is of partial Dini mean
oscillation with respect to $x'=(x_1,...,x_{n-1})$ in an open ball $B\subset\subset \Omega$,
if the function $\omega_f: (0, 1] \rightarrow [0,\infty)$ defined by
\[\omega_f(r)= \sup_{x\in B}\avint_{B_r(x)}
\left|f(y)-\avint_{B_{r}'(x')} f(z',y_n) dz'\right|dy,\quad B_{r}'(x')=\{y'\in \mathbb{R}^{n-1} : |y'-x'|<r\}\]
satisfies the Dini condition. 

\n iii) For each $\mathbf{f}\in L^1(\Omega)$, we define the following functions:
\[\varpi_{\mathbf{f}}(t)=\int_0^{t}\frac{\omega_{\mathbf{f}}(s)}{s}ds\text{ and }
\Phi_{\mathbf{f}}(t)=\varpi_{\mathbf{f}}(t)+\omega_{\mathbf{f}}(t)\]
\end{definition}

\begin{remark}\label{r1.1}

\n We observe that if a function $f\in L^1(\Omega)$ is such that for each $a_n$, the function 
$x' \rightarrow f(x',a_n)$ is H\"{o}lder continuous i.e. $|f(x',a_n)-f(y',a_n)|\leq C|x'-y'|^\alpha$ 
for some $\alpha\in(0,1)$, then it is easy to verify that $\omega_f(r)\leq 2Cr^\alpha$ for any $r\in(0,1]$,
which leads to
\[\int_0^1 \frac{\omega_{\mathbf{f}}(r)}{r}dr\leq \int_0^1 2Cr^{\alpha-1}dr=\frac{2C}{\alpha} <\infty\]
Hence, $f$ is of partial Dini mean oscillation with respect to $x'$ in any open ball $B\subset\subset \Omega$.
\end{remark}

\vs0.3cm \n Here is the main result of this paper:

\begin{theorem}\label{t1.1}
Assume that $\mathbf{A}$ and $\mathbf{f}$ satisfy
(1.1)-(1.3) and the following conditions:
\begin{align}\label{e1.5-1.6}
& \forall i,j=1,...,n, ~~a_{ij} \text{ is of partial Dini mean
oscillation with respect to } x' \text{ in } \Omega\\
& \forall i=1,...,n, ~~f_{i} \text{ is of partial Dini mean
oscillation with respect to } x' \text{ in } \Omega
\end{align}
Then for any weak solution of $(P)$, we have $u\in C^{0,1}_{loc}(\Omega)$.
\end{theorem}

\vs 0.3cm\n The novelty in Theorem 1.1 lies in the fact that Lipschitz continuity of weak solutions of problem $(P)$
is obtained even when the entries of the matrix $\mathbf{A}(x)$ and the vector function $\mathbf{f}(x)$
are discontinuous provided they satisfy a Dini mean oscillation condition in at least one direction
i.e. if they are regular in at least one variable. Since problem $(P)$ is invariant by rotation in the sense 
that it is transformed into a similar problem with different coefficients satisfying the same assumptions
as the original ones, it is obvious that we only need to have the Dini mean oscillation
condition satisfied in any arbitrary space direction.

\vs 0.3cm\n We recall that interior Lipschitz continuity for problem $(P)$ was established in 
\cite{[CL1]} and the same method was successfully extended to the quasilinear case in
\cite{[CL4]} and \cite{[CL5]}. Interior and boundary Lipschitz continuity were established in \cite{[L2]}
for a wide class of linear elliptic equations under some general assumptions.
Recently, in \cite{[L3]}, Lipschitz continuity was obtained using a different
method based on Harnack's inequality. This approach helped relax some of the assumptions required in
\cite{[CL1]} and \cite{[L2]} and only required that $\mathbf{A}(x) \in C^{0,\alpha}_{loc}(\Omega)$ and
$\text{div}(\mathbf{f}) \in L_{loc}^p(\Omega)$ for some $\alpha\in(0,1)$ and $p>n/(1-\alpha)$.

\vs 0.3cm\n Lastly, we would like to point out that the assumptions (1.4)-(1.5) were introduced in 
\cite{[DK]} to obtain $C^1$ and $C^2$-regularity of solutions to elliptic equations. In this regard, 
we also refer the reader to the recent work on gradient estimates for elliptic equations in 
divergence form with partial Dini mean oscillation coefficients \cite{[CKL]} .

\section{Estimates for the equation $\text{div}(\mathbf{A}(x)\nabla u)=-\text{div}(\mathbf{f})$}\label{2}

\n Under the assumptions of Theorem 1.1, it is known [see \cite{[C]}, Lemma 2.1]
that any weak solution $u$ of equation $\text{div}(\mathbf{A}(x)\nabla u)=-\text{div}(\mathbf{f})$ is such that
$~u\in C^{0,1}_{loc}(\Omega)$. The main result of this section is a local $L^\infty-$norm estimate of 
the gradient which will be used in the proof of Theorem 1.1 in section 3. Needless to say, this estimate 
is of interest for itself.

\begin{theorem}\label{t2.1} Let $\rho>0$ be such that $B_{3\rho}(x_0)\subset\subset\Omega$ and
$\Phi_{\mathbf{A}}(\rho)\leq C_0^*$ for some positive constant $C_0^*$ depending only on $n$ and $\lambda$.
Assume that $u\in H^1(B_{3\rho}(x_0))$ is a weak solution of equation
$\text{div}(\mathbf{A}(x)\nabla u)=-\text{div}(\mathbf{f})$ in $B_{3\rho}(x_0)$.
If moreover, we assume that $f_n\in L^\infty(B_{3\rho}(x_0))$
and $\mathbf{A}$ and $\mathbf{f}$ are of partial
Dini mean oscillation with respect to $x'$ in $B_{2\rho}(x_0)$, then 
$\nabla u\in  L^{\infty}(B_{\rho}(x_0))$ and we have for some positive constant $C_1$ depending only on $n$ and $\lambda$:
\[
|\nabla u|_{L^\infty(B_{\rho}(x_0))}\leq 3^{2n}C_1\left(\rho^{-n}|\nabla u|_{L^1(B_{3\rho}(x_0))}
+|f_n|_{L^\infty(B_{3\rho}(x_0))}+4\Phi_{\mathbf{f}}(\rho)\right)
\]
\end{theorem}

\n The proof of Theorem 2.1 requires a few lemmas.

\begin{lemma}\label{l2.1}
Assume that $\omega$ is a Dini function and let $a\in(0,1)$ and $b>1$ be two given real numbers.
Then the function defined by $\displaystyle{\widetilde{\omega}(t)=\sum_{i=0}^\infty a^i\left(\omega(b^it)\chi_{\{b^it\leq 1\}}+\omega(1)\chi_{\{b^it>1\}}\right) }$ satisfies
\[\int_0^t\frac{\widetilde{\omega}(s)}{s}ds\leq\frac{1}{1-a}\left[\varpi(t)+a\varpi(1)+\frac{\omega(1)}{\gamma }t^\gamma\right]
\quad\forall t\in[0,1]\]

\n where $\displaystyle{\gamma=-\frac{\ln(a)}{\ln(b)}}$ and $\displaystyle{\varpi(t)=\int_0^{t}\frac{\omega(s)}{s}ds}$. 
In particular, $\widetilde{\omega}$ is also a Dini function.
\end{lemma}

\vs 0,3cm\n \emph{Proof.} First, we recall that the function $\widetilde{\omega}$ was
introduced in [\cite{[D]}, Lemma 3.1], where an estimate was also given. Nevertheless, our estimate is new and more precise.

\n We start by writing $\widetilde{\omega}(t)=\widetilde{\omega}_1(t)+\widetilde{\omega}_2(t)$ for $t\in(0,1)$,
where
\begin{equation*}
\widetilde{\omega}_1(t)=\sum_{i=0}^\infty a^i\omega(b^it)\chi_{\{b^it\leq 1\}}\quad\text{and}\quad
\widetilde{\omega}_2(t)=\omega(1)\sum_{i=0}^\infty a^i\chi_{\{b^it>1\}}
\end{equation*}

\n Next, let $\displaystyle{i_0=\left[-\frac{\ln(t)}{\ln(b)}\right]}$ and observe that we have $t\leq b^{-i}$ iff $i\leq i_0$. 
Then we have
\begin{eqnarray}\label{e2.1}
\int_0^t \frac{\widetilde{\omega}_1(s)}{s}ds&=& \int_0^t a^0\frac{\omega(b^0s)}{s}\chi_{\{b^0s\leq 1\}}ds+\sum_{i=1}^\infty a^i \int_0^t\frac{\omega(b^is)}{s}\chi_{\{b^is\leq 1\}}ds\nonumber\\
&=&\int_0^t\frac{\omega(s)}{s}ds+\sum_{i=1}^\infty a^i \int_0^{tb^{i}}\frac{\omega(\tau)}{\tau}\chi_{\{\tau\leq 1\}}d\tau\nonumber \\
&\leq& \varpi(t)+\sum_{i=1}^\infty a^i \int_0^{1}\frac{\omega(\tau)}{\tau}d\tau=\varpi(t)+\frac{a}{1-a}.\varpi(1)
\end{eqnarray}

\n Given that $\displaystyle{i_0\leq-\frac{\ln(t)}{\ln(b)}<i_0+1}$, we can write
\begin{eqnarray*}
\widetilde{\omega}_2(s)= \omega(1)\sum_{i=0}^\infty a^i\chi_{\{i>-\frac{\ln(t)}{\ln(b)}\}}
= \omega(1)\sum_{i=i_0+1}^\infty a^i =\frac{a^{i_0+1}}{1-a}\omega(1)
\end{eqnarray*}

\n Moreover, since $a<1$, we have
$\displaystyle{a^{i_0+1}\leq a^{-\frac{\ln(t)}{\ln(b)}}=t^\gamma}$, which leads to
$\displaystyle{\widetilde{\omega}_2(s)\leq \frac{\omega(1)}{1-a}t^{\gamma}}$ and

\begin{eqnarray}\label{e2.2}
\int_0^t \frac{\widetilde{\omega}_2(s)}{s}ds&\leq& \frac{\omega(1)}{1-a}\int_0^t\frac{s^{\gamma}}{s}ds~\leq~\frac{\omega(1)}{\gamma(1-a)}t^\gamma
\end{eqnarray}

\n Now, combining (2.1) and (2.2), we obtain
\begin{eqnarray*}
\int_0^t \frac{\widetilde{\omega}(s)}{s}ds&\leq& \varpi(t)+\frac{a}{1-a}.\varpi(1)+\frac{\omega(1)}{\gamma(1-a)}t^\gamma\\
&\leq&\frac{1}{1-a}\left[(1-a)\varpi(t)+a\varpi(1)+\frac{\omega(1)}{\gamma }t^\gamma\right]\\
&\leq&\frac{1}{1-a}\left[\varpi(t)+a\varpi(1)+\frac{\omega(1)}{\gamma }t^\gamma\right]
\end{eqnarray*}
\qed

\n The following lemma is a slight improvement of Theorem 1.2 of \cite{[CKL]}
in the sense that it provides a more precise $L^\infty$-estimate of the gradient.

\begin{lemma}\label{l2.2}
Let $u\in H^1(B_3)$ be a weak solution of the equation
$\text{div}(\mathbf{A}(x)\nabla u)=-\text{div}(\mathbf{f})$ in $B_3$,
with $\mathbf{A}$ and $\mathbf{f}$ satisfying (1.1)-(1.3) in $B_3$
and both $\mathbf{A}$ and $\mathbf{f}$ of partial
Dini mean oscillation with respect to $x'$ in $B_2$. Then we have
$\nabla u\in  L^{\infty}(B_1)$ with

\[|\nabla u|_{L^\infty(B_1)}\leq 3^{nk_0}C_1\left(|\nabla u|_{L^1(B_3)}
+\overline{f}+4\Phi_{\mathbf{f}}(1)\right)\]
where $C_1=C_1(n,\lambda)$ is a positive constant depending only on $n$ and $\lambda$,
and $k_0$ is an integer greater than 1 satisfying
\[\varpi_{\mathbf{A}}(2^{-k_0})+2\omega_{\mathbf{A}}(1)\sqrt{2^{-k_0}}\leq 3^{-n-1}C_0^{-1}=C_0^*\]
\end{lemma}

\n \emph{Proof.} First, we observe that by scaling, we may replace the ball $B_3$ by $B_6$
as in \cite{[CKL]}. Next, we denote by $C_0=C_0(n,\lambda)$ the positive constant depending only on 
$n$ and $\lambda$ that was introduced in [\cite{[CKL]}, Proof of Theorem 1.2, p. 1520]. 
Following this reference, we choose $\displaystyle{\gamma=\frac{1}{2}}$, $0<\kappa<\min(2^{-1},C_0^{-2})$ and we denote by
$\widetilde{\omega}_{\mathbf{A}}(t)$ the function defined in Lemma 2.1 with $a=\sqrt{\kappa}$
and $\displaystyle{b=\frac{1}{\kappa}}$. Let now $k_0$ be a positive integer greater than 1 that satisfies
\[C_0\int_0^{\frac{1}{2^{k_0}}}\frac{\widetilde{\omega}_{\mathbf{A}}(t)}{t}dt\leq 3^{-n}\]

\n which by taking into account the estimate of Lemma 2.1 is true if

\[\frac{C_0}{1-\sqrt{\kappa}}\left(\varpi_{\mathbf{A}}(2^{-k_0})+\sqrt{\kappa}\varpi_{\mathbf{A}}(1)+2\omega_{\mathbf{A}}(1)\sqrt{2^{-k_0}}\right)\leq 3^{-n}\]

\n At this step, we further assume that $\kappa\leq 2^{-2}$, which leads to $\displaystyle{1-\sqrt{\kappa}\geq 2^{-1} }$, and makes the above inequality hold if
\[2(\varpi_{\mathbf{A}}(2^{-k_0})+\sqrt{\kappa}\varpi_{\mathbf{A}}(1)+2\omega_{\mathbf{A}}(1)\sqrt{2^{-k_0}})\leq
3^{-n}C_0^{-1}\]

\n This in turn remains true if $\kappa$ and $k_0$ are chosen such that 
\[2\sqrt{\kappa}\varpi_{\mathbf{A}}(1),~ \varpi_{\mathbf{A}}(2^{-k_0})+2\omega_{\mathbf{A}}(1)\sqrt{2^{-k_0}}\leq 3^{-n-1}C_0^{-1}=C_0^*\]

\n If we replace $f_1$ by $f_n$, we get the estimate [see \cite{[CKL]}, p. 1522] with a positive constant $C_1(n,\lambda)$ depending only on $n$ and $\lambda$
\begin{eqnarray*}
3^{-nk_0}|\nabla u|_{L^\infty(B_1)}&\leq& C_1(n,\lambda)\left(|\nabla u|_{L^1(B_3)}+|f_n|_{L^\infty(B_3)}+\int_0^1 \frac{\widetilde{\omega}_{\mathbf{f}}(t)}{t}dt\right)
\end{eqnarray*}

\n which can be written by using Lemma 2.1 again as
\begin{eqnarray*}
|\nabla u|_{L^\infty(B_1)}&\leq& 
3^{nk_0}C_1\left(|\nabla u|_{L^1(B_3)}+\overline{f}+\frac{1}{1-\sqrt{\kappa}}(\varpi(1)+a\varpi(1)+2\omega(1)1^\gamma)\right)\nonumber\\
&\leq& 3^{nk_0}C_1\left(|\nabla u|_{L^1(B_3)}+\overline{f}+4(\omega_{\mathbf{f}}(1)+\varpi_{\mathbf{f}}(1))\right)\nonumber\\
&=& 3^{nk_0}C_1\left(|\nabla u|_{L^1(B_3)}+\overline{f}+4\Phi_{\mathbf{f}}(1)\right)
\end{eqnarray*}
\qed

\n The following lemma is a slight improvement of Lemma 2.2 of \cite{[L3]}.

\begin{lemma}\label{l2.3} Assume that $u$ is a nonnegative weak solution of the equation
$\text{div}(\mathbf{A}(x)\nabla u)=-\text{div}(\mathbf{f})$ in $\Omega$ and let $x_0\in \Omega$ and $r>0$ such that
$B_{5r}(x_0)\subset\subset \Omega$ and $\overline{B}_r(x_0)\cap\{u=0\}\neq \emptyset$. Then we
have for some positive constant $C_2$ depending only on $n$, $\lambda$
and $\bar f$: $\qquad\displaystyle{\max_{\overline{ B_{r}}(x_0)} u \,\leq \,C_2\,r}.$
\end{lemma}

\n \emph{Proof.} Let $x_1\in \overline{B}_r(x_0)\cap\{u=0\}$ and let $\omega_n=|B_1|$ 
be the measure of the unit ball in $\mathbb{R}^n$.
First since $B_{5r}(x_0)\subset\subset \Omega$, it is easy to verify that
$B_{4r}(x_1)\subset\subset \Omega$.
Next, we observe that since $\mathbf{f}\in L^{\infty}(\Omega)$, we can apply Harnack's
inequality [\cite{[GT]}, Theorem 8.17-Theorem 8.18, p. 194] to equation
$\text{div}(\mathbf{A}(x)\nabla u)=-\text{div}(\mathbf{f})$  with $p=n+1$. Therefore, 
we get for a positive constant $C$ depending only on $n$
\begin{eqnarray*}
\max_{\overline{B}_{2r}(x_1)} ~ u&\leq& C
\Big(\min_{\overline{B}_{2r}(x_1)} ~ u +
{1\over\lambda}r^{1-{n\over{n+1}}}||\mathbf{f}||_{n+1,\overline{B}_{2r}(x_1)}\Big)\nonumber\\
&\leq& C\Big(0+{\overline{f}\over\lambda}r^{{1\over {n+1}}}.|B_{2r}(x_1)|^{1\over{n+1}}\Big)\nonumber\\
&=& \frac{C\overline{f}2^{{n\over{n+1}}}\omega_n^{1\over{n+1}}}{\lambda}.r^{{1\over{n+1}}}.r^{{n\over{n+1}}}
=\frac{C\overline{f}2^{{n\over{n+1}}}\omega_n^{1\over{n+1}}}{\lambda}r=C_2 r
\end{eqnarray*}

\n Given that $\overline{B}_r(x_0)\subset \overline{B}_{2r}(x_1)$, the lemma follows.
\qed

\vs 0,3cm\n The following lemma is a Cacciopoli type lemma.

\begin{lemma}\label{l2.4} Assume that $u$ is a weak solution of the equation
$\text{div}(\mathbf{A}(x)\nabla u)=-\text{div}(\mathbf{f})$ in $\Omega$.
For each open ball $B_r(x_0)$ such that $B_{2r}(x_0)\subset \Omega$, we have:
\[\int_{B_{r}(x_0)}|\nabla u|^2dx~\leq~\frac{32}{\lambda^4r^2}\int_{B_{2r}(x_0)}u^2 dx
+\frac{8\overline{h}}{\lambda r}\int_{B_{2r}(x_0)}  |u| dx
+\frac{2^{n+1}\overline{f}^2\omega_n}{\lambda^2} r^n\]
\end{lemma}

\n \emph{Proof.} Let $B_r(x_0)$ be an open ball such that $B_{2r}(x_0)\subset \Omega$,
and let $\eta\in C_0^\infty(B_{2r}(x_0))$ be a cut-off function such that
$$\eta=1 \quad\hbox{in}\quad B_{r}(x_0), \qquad 0\leq \eta\leq 1\qquad\hbox{and}\qquad |\nabla \eta|\leq \frac{2}{r} \quad\hbox{in}\quad B_{2r}(x_0).$$
\n Using $\eta^2u$ as a test function for equation $(P)ii)$, we get
\begin{eqnarray*}
\int_{B_{2r}(x_0)}a(x) \nabla u.\nabla(\eta^2u) dx=-\int_{B_{2r}(x_0)}\mathbf{f}(x).\nabla(\eta^2 u) dx
\end{eqnarray*}
which can be written as
\begin{eqnarray*}
&&\int_{B_{2r}(x_0)}\eta^2a(x) \nabla u.\nabla udx~=~-\int_{B_{2r}(x_0)}2\eta.u. a(x)\nabla u.\nabla\eta dx\\
&&\quad-\int_{B_{2r}(x_0)} 2\eta u.\mathbf{f}(x).\nabla \eta dx-\int_{B_{2r}(x_0)}\eta^2.\mathbf{f}(x) .\nabla u dx
\end{eqnarray*}

\n or by using (1.1)-(1.2)
\begin{eqnarray*}
&&\lambda\int_{B_{2r}(x_0)}\eta^2|\nabla u|^2dx~\leq~\frac{2}{\lambda}\int_{B_{2r}(x_0)}\eta.|\nabla u|.|u|.|\nabla\eta| dx\\
&&\quad+\int_{B_{2r}(x_0)} 2\eta |u|.|\mathbf{f}|.|\nabla \eta| dx+\int_{B_{2r}(x_0)}\eta.|\nabla u|.\eta.|\mathbf{f}|_\infty dx
\end{eqnarray*}
\n By taking into account (1.3) and the fact that $|\nabla \eta|\leq 2/r$, we obtain
\begin{eqnarray*}
&&\int_{B_{2r}(x_0)}\eta^2|\nabla u|^2dx~\leq~\int_{B_{2r}(x_0)}(\eta|\nabla u|).\left(\frac{4|u|}{\lambda^2r}\right) dx\\
&&\quad+\frac{4\overline{f}}{\lambda r}\int_{B_{2r}(x_0)}  |u| dx+\int_{B_{2r}(x_0)}(\eta|\nabla u|).\left(\frac{\overline{f}}{\lambda}\right) dx
\end{eqnarray*}
\n Now, we apply the following Young's type inequality $\displaystyle{ab\leq \frac{1}{4}a^2+b^2}$
to the first and third integrals of the righthand side of the previous inequality
\begin{eqnarray*}
&&\int_{B_{2r}(x_0)}\eta^2|\nabla u|^2dx~\leq~\frac{1}{4}\int_{B_{2r}(x_0)}\eta^2|\nabla u|^2dx+
\frac{16}{\lambda^4r^2}\int_{B_{2r}(x_0)}u^2 dx\\
&&\quad+\frac{4\overline{h}}{\lambda r}\int_{B_{2r}(x_0)}  |u| dx+
\frac{1}{4}\int_{B_{2r}(x_0)}\eta^2|\nabla u|^2dx
+\int_{B_{2r}(x_0)}\frac{\overline{f}^2}{\lambda^2} dx
\end{eqnarray*}
\n which leads, since $\eta=1$ in $B_{r}(x_0)$, to

\begin{eqnarray*}
&&\int_{B_{r}(x_0)}|\nabla u|^2dx~\leq~\frac{32}{\lambda^4r^2}\int_{B_{2r}(x_0)}u^2 dx
+\frac{8\overline{f}}{\lambda r}\int_{B_{2r}(x_0)} |u| dx
+\frac{2^{n+1}\omega_n\overline{f}^2}{\lambda^2} r^n
\end{eqnarray*}
\qed

\vs 0,3cm\n Combining Lemmas 2.3 and 2.4, we obtain the following lemma.

\begin{lemma}\label{l2.5} Assume that $u$ is a nonnegative weak solution of the equation
$\text{div}(\mathbf{A}(x)\nabla u)=-\text{div}(\mathbf{f})$ in $\Omega$ and let $x_0\in \Omega$ and $r>0$ such that
$B_{10r}(x_0)\subset\subset \Omega$ and $\overline{B}_{2r}(x_0)\cap\{u=0\}\neq \emptyset$. Then we
have:
\[\int_{B_{r}(x_0)}|\nabla u|dx~\leq~C_3r^n\]
where $\displaystyle{C_3=\frac{\omega_n2^{\frac{n+1}{2}}}{\lambda^2}\sqrt{16C_2^2
+4C_2\overline{f}\lambda^3+\overline{f}^2\lambda^2}}$ and $C_2$ is the constant in Lemma 2.3.
\end{lemma}

\n \emph{Proof.} From Lemmas 2.3 and 2.4, we have
\begin{eqnarray}\label{e2.3-2.4}
&&\max_{\overline{ B_{2r}}(x_0)} u \,\leq \,C_2(2r)\text{ and }\\
&&\int_{B_{r}(x_0)}|\nabla u|^2dx~\leq~\frac{32}{\lambda^4r^2}\int_{B_{2r}(x_0)}u^2 dx
+\frac{8\overline{f}}{\lambda r}\int_{B_{2r}(x_0)}  u dx
+\frac{2^{n+1}\overline{f}^2\omega_n}{\lambda^2} r^n
\end{eqnarray}

\n Combining (2.3) and (2,4), we obtain
\begin{eqnarray}\label{e2.5}
\int_{B_{r}(x_0)}|\nabla u|^2dx&\leq&\frac{32}{\lambda^4 r^2}\int_{B_{2r}(x_0)}C_2^2(2r)^2 dx
+\frac{8\overline{f}}{\lambda r}\int_{B_{2r}(x_0)}  C_2(2r) dx
+\frac{2^{n+1}\omega_n\overline{f}^2}{\lambda^2} r^n\nonumber\\
&\leq&\frac{128}{\lambda^4}\omega_n C_2^2(2r)^n
+\frac{16\overline{f}}{\lambda}\omega_nC_2(2r)^n
+\frac{2^{n+1}\overline{f}^2\omega_n}{\lambda^2} r^n\nonumber\\
&=&\frac{2^{n+1}\omega_n}{\lambda^4}\left[64C_2^2
+8C_2\overline{f}\lambda^3+\overline{f}^2\lambda^2\right]r^n
\end{eqnarray}

\n We conclude by using Cauchy-Schwartz inequality and taking into account (2.5)
\begin{eqnarray*}
\int_{B_{r}(x_0)}|\nabla u|dx&\leq&\left(\int_{B_{r}(x_0)}|\nabla u|^2dx\right)^{\frac{1}{2}}.|B_{r}(x_0)|^{\frac{1}{2}}
\leq C_3r^n
\end{eqnarray*}
\qed

\n \emph{Proof of Theorem 2.1.} We observe that the function $\displaystyle{v(y)={u(x_0+\rho y)\over \rho}}$
satisfies the equation $\text{div}(\mathbf{A}_\rho(x)\nabla v)=-\text{div}(\mathbf{f}_\rho)$ in $B_{3}$,
where $\mathbf{F}_\rho(y)=\mathbf{F}(x_0+\rho y)$. Moreover, it is obvious that $\mathbf{A}_\rho$ and
$\mathbf{f}_\rho$ satisfy the assumption of Lemma 2.2 in $B_3$. Therefore, we get the estimate
\begin{equation}\label{e2.6}
|\nabla v|_{L^\infty(B_1)}\leq 3^{nk_0}C_1\left(|\nabla v|_{L^1(B_3)}+|f_{\rho n}|_{L^\infty(B_3)}+
4\Phi_{\mathbf{f}_\rho}(1)\right)
\end{equation}

\n where $C_1$ is a positive constant depending only on $n$ and $\lambda$,
$k_0$ is an integer satisfying
\[\varpi_{\mathbf{A}_\rho}(2^{-k_0})+2\omega_{\mathbf{A}_\rho}(1)\sqrt{2^{-k_0}}\leq C_0^*\]

\n and $C_0$ is the positive constant depending only on $n$ and $\lambda$ 
from [\cite{[CKL]}, p. 1520].

\n If we observe that $\omega_{\mathbf{f}_\rho}(s)=\omega_{\mathbf{f}}(\rho s)$, 
$\varpi_{\mathbf{f}_\rho}(s)=\varpi_{\mathbf{f}}(\rho s)$, and we choose $k_0=2$, this condition reduces to  
$\varpi_{\mathbf{A}}(2^{-2}\rho)+\omega_{\mathbf{A}}(\rho)\leq C_0^*$, which is
in particular true if $\Phi_{\mathbf{A}}(\rho)\leq C_0^*$.

\n Finally, since $\Phi_{\mathbf{f}_\rho}(s)=\Phi_{\mathbf{f}}(\rho s)$, the estimate of Theorem 2.1 follows from (2.6).
\qed

\section{Proof of Theorem 1.1}\label{2}

\n Let $\epsilon>0$, $\Omega_\epsilon=\{x\in \Omega\,/\,d(x,\partial\Omega)>\epsilon\,\}$, and let
$C_0^*$ be the positive constant depending only on $n$ and $\lambda$ that was introduced in Theorem 2.1. Since
$\displaystyle{\lim_{t\rightarrow0} \Phi_{\mathbf{F}}(t)=0}$, there exists $t_0\in(0,1)$ such that
\begin{equation}\label{e3.1}
\Phi_{\mathbf{A}}(t)\leq C_0^*\quad\forall t\in(0,t_0)
\end{equation}

\n We fix $\displaystyle{\epsilon_0=\frac{3}{4}t_0}$ and assume that $\epsilon<\epsilon_0$. 
We shall prove that $\nabla u$ is bounded in $\Omega_{41\epsilon}$ by a constant depending only on $n$, 
$\lambda$, $\overline{f}$, $\Phi_{\mathbf{f}}(1)$, $|\nabla u|_{L^1(\Omega)}$ and $\epsilon$, provided 
$\epsilon<\epsilon_0$.

\n Let $x_0\in \Omega_{41\epsilon}$. We distinguish two cases:

\vs0.2cm \n i)
\underline{$B_{3\epsilon}(x_0)\subset\{u>0\}$}:

\vs0.2cm \n Since $\chi=1$ a.e. in $B_{3\epsilon}(x_0)$, $u$ satisfies the equation
$\text{div}(\mathbf{A}(x)\nabla u)=-\text{div}(\mathbf{f})$ in $B_{3\epsilon}(x_0)$.
We also have by (1.3)-(1,5) that $||f_n||_{\infty}\leq \overline{f}$
and $\mathbf{A}$ and $\mathbf{f}$ are of partial
Dini mean oscillation with respect to $x'$ in $B_{2\epsilon}(x_0)$.
Moreover, we have from (3.1) $\Phi_{\mathbf{A}}(\epsilon)\leq C_0^*$.
Therefore, by Theorem 2.1 applied with $\rho=\epsilon$, we get for a positive constant $C_1=C_1(n,\lambda)$
\begin{eqnarray}\label{e3.2}
|\nabla u|_{L^\infty(B_{\epsilon}(x_0))}
&\leq& 3^{2n}C_1\left(\epsilon^{-n}|\nabla u|_{L^1(\Omega)}+|f_n|_{L^\infty(B_{3\epsilon}(x_0))}+4\Phi_{\mathbf{f}}(\epsilon)\right)\nonumber
\end{eqnarray}

\n Since $\Phi_{\mathbf{f}}(t)$ is nondecreasing and $\epsilon<\epsilon_0<1$, we obtain
\begin{eqnarray}\label{e3.2}
|\nabla u|_{L^\infty(B_{\epsilon}(x_0))}
&\leq& 3^{2n}C_1\left(\epsilon^{-n}|\nabla u|_{L^1(\Omega)}+\overline{f}+4\Phi_{\mathbf{f}}(1)\right)
\end{eqnarray}

\vs0.2cm \n ii)
\underline{$B_{3\epsilon}(x_0)\cap\{u=0\}\neq\emptyset$}:

\vs0.2cm \n Let $x\in B_\epsilon(x_0)$ such that
$u(x)>0$ and let $r(x)=\text{dist}(x,\{u=0\})$ be the distance function to
the set $\{u=0\}$. Our objective is to estimate $|\nabla u|_{L^\infty(B_{r(x)/3}(x))}$.
To do that, we will again use Theorem 2.1.

\n We claim that $r(x)<4\epsilon$ and $\overline{B_{10r(x)}}(x) \subset B_{41\epsilon}(x_0)$.
Indeed, by assumption ii), there exists $z\in
B_{3\epsilon}(x_0)\cap \{u=0\}$. So, we get $ r(x)\leq  |x-z| \leq |x-x_0|+ |x_0-z|< \epsilon+3\epsilon=4\epsilon$.

\n Then we have for each $y\in \overline{B_{10r(x)}}(x)$
$$|x_0-y|\leq |x_0-x|+|x-y|< \epsilon+10r(x)<\epsilon+10(4\epsilon)=41\epsilon$$

\n This means that $y\in B_{41\epsilon}(x_0)$, and therefore $\overline{B_{10r(x)}}(x) \subset B_{41\epsilon}(x_0)$.
In particular, we have $\overline{B_{10r(x)}}(x) \subset \Omega_{41\epsilon}\subset \Omega$.

\vs0.1cm\n Now, we obtain from the pevious step that
$\overline{B_{r(x)}}(x)\subset\Omega$. 
Moreover, since $B_{r(x)}(x) \subset\{u>0\}$, we have by $(P)i)$ that $\chi=1$ a.e. in $B_{r(x)}(x)$,
which leads by $(P)ii)$ to $\text{div}(\mathbf{A}(x)\nabla u)=-\text{div}(\mathbf{f})$ in $B_{r(x)}(x)$.
By (1.3)-(1,5), we also know that $||f_n||_{\infty}\leq \overline{f}$
and $\mathbf{A}$ and $\mathbf{f}$ are of partial
Dini mean oscillation with respect to $x'$ in $B_{2r(x)/3}(x)$.
On the other hand, we have $r(x)/3\leq 4\epsilon/3<t_0$,
which ensures by (3.1) that $\Phi_{\mathbf{A}}(r(x)/3)\leq C_0^*$.
Hence, we infer from Theorem 2.1 applied with $\rho=r(x)/3$, that we have for a positive constant $C_1=C_1(n,\lambda)$

\begin{eqnarray}\label{e3.3}
|\nabla u|_{L^\infty(B_{r(x)/3}(x))}
&\leq& 3^{2n}C_1\left(\left(\frac{r(x)}{3}\right)^{-n}|\nabla u|_{L^1(B_{r(x)})}+\overline{f}+4\Phi_{\mathbf{f}}(r(x)/3)\right)
\end{eqnarray}

\n Since $r(x)/3<t_0<1$ and $\Phi_{\mathbf{f}}(t)$ is nondecreasing, we have
$\Phi_{\mathbf{f}}(r(x)/3)\leq \Phi_{\mathbf{f}}(1)$. Furthermore, we have $\overline{B_{10r(x)}}(x)\subset \Omega$ 
and $\overline{B}_{2r(x)}(x) \cap\{u=0\}\neq\emptyset$, whence we can use Lemma 2.5 to improve (3.3) as follows
\begin{eqnarray*}
|\nabla u|_{L^\infty(B_{r(x)/3}(x))}~\leq~ 3^{2n}C_1\left(3^{n}C_3+\overline{f}+4\Phi_{\mathbf{f}}(1)\right)
\end{eqnarray*}

\vs0.2cm \n Given that $x$ is arbitrary in $B_\epsilon(x_0)\cap\{u>0\}$ and $\nabla u(x)=0$ a.e. in $B_\epsilon(x_0)\cap\{u=0\}$, 
it turns out that $\nabla u$ is uniformly bounded in $B_\epsilon(x_0)$, with
\begin{eqnarray}\label{e3.4}
|\nabla u|_{L^\infty(B_{\epsilon}(x_0))}&\leq&
 3^{2n}C_1\left(3^{n}C_3+\overline{f}+4\Phi_{\mathbf{f}}(1)\right)
\end{eqnarray}

\n Finally, by taking into account the fact that $\displaystyle{ \Omega_{41\epsilon}\subset\bigsqcup_{0<\epsilon<\epsilon_0,x_0\in \Omega_{41\epsilon}} B_\epsilon(x_0) }$, we conclude
from (3.2) and (3.4) that 
\begin{eqnarray*}
|\nabla u|_{L^\infty(\Omega_{41\epsilon})}&\leq&
 3^{2n}C_1\left(\epsilon^{-n}|\nabla u|_{L^1(\Omega)}+3^{n}C_3+\overline{f}+4\Phi_{\mathbf{f}}(1)\right)
\end{eqnarray*}

\n which means that $|\nabla u(x)|$ is uniformly bounded in $\Omega_{41\epsilon}$ by a constant depending only on $n$, $\lambda$,
$\overline{f}$, $\Phi_{\mathbf{f}}(1)$, $|\nabla u|_{L^1(\Omega)}$ and $\epsilon<d(x,\partial\Omega)/41$, for any  
small enough $\epsilon>0$.
\qed

\begin{remark}\label{r3.1}

With a slight modification of the proof, it is not difficult to extend 
Theorem 1.1 to the following problem:
\begin{equation*}
\begin{cases}
& \text{Find } (u, \chi) \in  H^{1}(\Omega)\times L^\infty (\Omega)
\text{ such that}:\\
& (i)\quad  u\geq 0, \quad 0\leq \chi\leq 1 ,
\quad u(1-\chi) = 0 \,\,\text{ a.e.  in } \Omega\\
& (ii)\quad \text{div}(\mathbf{A}(x) \nabla
u + \chi \mathbf{f}(x))=-\chi \mathbf{g}(x)\quad \text{ in } \Omega
\end{cases}
\end{equation*}

\n provided that $\mathbf{A}(x)$ satisfies (1.1)-(1.2), $\mathbf{g}(x)$ and $\mathbf{f}(x)$ 
satisfy (1.3), and the three functions are of partial Dini mean
oscillation with respect to $x'$ in $\Omega$.
\end{remark}

\end{document}